\documentclass[a4paper,12pt,twoside]{article}
\usepackage{times}
\usepackage[english]{babel}
\usepackage{amssymb,amsmath}
\numberwithin{equation}{section}
\usepackage{amsthm}
\usepackage{array}
\usepackage[noadjust]{cite}
\usepackage{hyperref}
\usepackage[height=22.7cm , width = 16cm , top = 4cm , left = 3cm, a4paper]{geometry}
\usepackage{amssymb}
\usepackage{amsmath}
\usepackage{cite}
 \usepackage{pstricks}
\usepackage{epsfig}
\usepackage{verbatim}
\usepackage{multicol}
\usepackage{graphicx, color, psfrag}
\usepackage{graphics, psfrag}
\usepackage{color,soul}
\newtheorem{theorem}{Theorem}[section]
\newtheorem{remark}[theorem]{Remark}
\newcommand{\e}{\end{document}}

\newtheorem{corollary}{Corollary}[section]

\theoremstyle{definition}
\newtheorem{definition}{Definition}[section]
\begin{document}

\thispagestyle{empty}

\author{F. M. Abdel Moneim$^{(1)}$\thanks{Corresponding author, email: \href{mailto:fatmakasem1982@gmail.com}{fatmakasem1982@gmail.com}}~, Abdelfattah Mustafa $^{(1,2)}$\thanks{ Email: \href{mailto: amelsayed@mans.edu.eg}{amelsayed@mans.edu.eg}} and B. S. El-Desouky$^{(1)}$ \\
$^{(1)}$ {\small Department of Mathematics, Faculty of Science,
Mansoura University, Mansoura 35516, Egypt.} \\
$^{(2)}$ {\small Department of Mathematics, Faculty of Science, Islamic University of Madinah, KSA.}}

\title{New Generalization Families of Higher Order Changhee Numbers and Polynomials}

\date{}

\maketitle
\small \pagestyle{myheadings}
        \markboth{{\scriptsize New Generalization Families of Higher Order Changhee Numbers and Polynomials}}
        {{\scriptsize {F. M. Abdel Moneim, A. Mustafa and B. S. El-Desouky}}}

\hrule \vskip 8pt

\begin{abstract}
In this paper, we introduce new generalization of higher order Changhee of the first and second kind. Moreover, we derive some new results for these numbers and polynomials. Furthermore, some interesting special cases of the generalized higher order Changhee numbers and polynomials are deduced.
\end{abstract}



\section{Introduction}
The usual Stirling numbers and the singles Stirling numbers of the first kind, $s(n,k)$  and $s_1 (n,k)$ are defined, respectively, by
\begin{eqnarray} \label{1.1}
(x)_n &=& \sum_{k=0}^{n} s(n,k) x^k, \quad s(n,0)=\delta_{n,0}, \; \text{and} \; s(n,k)=0, \; \text{for} \; k>n.
\\ \label{eq:1.2}
\langle x \rangle_n &=& \sum_{k=0}^{n} s_1(n,k) x^k, \quad s_1(n,0)=\delta_{n,0}, \; \text{and} \; s_1(n,k)=0, \; \text{for} \; k>n,
\end{eqnarray}
where  $(x)_n = x(x-1)\cdots (x-n+1)$  , and   $\langle x \rangle_n=x(x+1)\cdots (x+n-1)$.\\

\noindent
These numbers satisfy the following recurrence relations
\begin{equation}\label{1.3}
s(n+1,k)=s(n,k-1)-ns(n,k).
\end{equation}

\noindent
The generalized Comtet numbers of the first and second kind, $s_{\bar{\pmb \alpha}}  (n,i;\bar{\pmb r} )$  and  $S_{\bar{\pmb \alpha}}  (n,i;\bar{\pmb r} )$, see \cite{2}, are defined, respectively, by
\begin{equation} \label{1.4}
(x;\bar{\pmb \alpha},\bar{\pmb r} )_n= \sum_{i=0}^n  s_{\bar{\pmb \alpha}}  (n,i; \bar{\pmb r})   x^i,
\end{equation}

\noindent
and
\begin{equation} \label{1.5}
x^n= \sum_{i=0}^n   S_{\bar{\pmb \alpha}}  (n,i;\bar{\pmb r} )(x; \bar{\pmb \alpha}, \bar{\pmb r})_i ,
\end{equation}

\noindent
where $(x;\bar{\pmb \alpha}, \bar{\pmb r})_n=\prod_{i=0}^{n-1} (x-\alpha_i  )^{r_i }, \bar{\pmb \alpha }=(\alpha_0,\alpha_1,\cdots,\alpha_{n-1}), \bar{\pmb r}=(r_0,r_1,\cdots,r_{n-1} )$ .

\noindent
The $n$-th Changhee polynomials are defined by the generating function, \cite{5} -- \cite{12},
\begin{equation} \label{1.6}
\frac{2}{t+2}  (1+t)^x= \sum_{n=0}^\infty  Ch_n  (x)  \frac{t^n}{n!}.
\end{equation}

\noindent
In the special case, $x=0$ ,$Ch_n=Ch_n (0)$ are called Changhee numbers, and for $n\geq 0$,
\begin{equation} \label{1.7}
\int_{\mathbb{Z}_p} (x )_n  d\mu_{-1} (x)=Ch_n.
\end{equation}

\noindent
For $k\in \mathbb{N}$, Kim \cite{5} introduced Changhee numbers of the first kind of order $k$ by
\begin{equation} \label{1.8}
Ch_n^{(k)}=\int_{\mathbb{Z}_p} \cdots \int_{\mathbb{Z}_p}(x_1+x_2+\cdots+x_k)_n  d\mu_{-1} (x_1 )\cdots d\mu_{-1} (x_k ),
\end{equation}

\noindent
where $n$ is nonnegative integer.\\

\noindent
The generating function of these numbers are given as
\begin{equation} \label{1.9}
\sum_{n=0}^\infty Ch_n^{(k)}    \frac{t^n}{n!}=\left(\frac{2}{2+t}\right)^k,
\end{equation}

\noindent
where $n\in \mathbb{Z}\geq 0, k\in \mathbb{N}$.

\noindent
The higher-order Changhee polynomials are defined by
\begin{equation} \label{1.10}
Ch_n^{(k)} (x)=\int_{\mathbb{Z}_p}\cdots \int_{\mathbb{Z}_p} (x_1+x_2+\cdots +x_k+x)_n  d\mu_{-1} (x_1 )\cdots d\mu_{-1} (x_k ).
\end{equation}

\noindent
For $k\in \mathbb{Z}$, the Euler polynomials of order $k$ are defined by the generating function to be, \cite{1,5,12},
\begin{equation} \label{1.11}
\left(\frac{2}{e^t+1} \right)^k e^{xt}= \sum_{n=0}^\infty E_n^{(k) } (x)   \frac{t^n}{n!},
\end{equation}

\noindent
where  $x=0, E_n^{(k)}=E_n^{(k) } (0)$ are called the Euler numbers of order $k$ .

\noindent
Also, Kim \cite{7}, proved that
\begin{equation}\label{1.12}
Ch_n^{(k)} (x)= \sum_{\ell=0}^n s(n,\ell)    E_\ell^{(k)} (x),
\end{equation}
and
\begin{equation} \label{1.13}
E_n^{(k)} (x)= \sum_{\ell=0}^n S(n,\ell)    Ch_\ell^{(k)} (x) .
\end{equation}

\noindent
An explicit formula for higher order Changhee numbers are given by
\begin{equation} \label{1.14}
Ch_n^{(k)}=\left(-\frac{1}{2} \right)^n \sum_{\ell=0}^n s(n,\ell)  (k+n-1)^\ell ,
\end{equation}

\noindent
where $s(n,\ell)$ are Stirling numbers of the first kind, see \cite{2,5,7,8}.

\section{Multiparameter Changhee Numbers of the First Kind}
In this section, the new definitions for Changhee numbers and polynomials are introduced. Some new results are derived; also some special cases are established as follows.
\begin{definition}\label{def: 2.1}
The multiparameter Changhee numbers   ${\check{C}}h_{n;\bar{\pmb \alpha},\bar{\pmb r}}^{(k)}$    are defined by
\begin{equation} \label{2.1}
\check{C}h_{n;\bar{\pmb \alpha}, \bar{\pmb r}}^{(k) }=\int_{\mathbb{Z}_p} \cdots \int_{\mathbb{Z}_p} \prod_{i=0}^{n-1} (x_1 x_2\cdots x_k-\alpha_i  )^{r_i}    d\mu_{-1} (x_1 )\cdots d\mu_{-1} (x_k ),
\end{equation}

\noindent
where $\bar{\pmb \alpha}= (\alpha_0,\alpha_1,\cdots ,\alpha_{n-1} ), \bar{\pmb r} = (r_0,r_1,\cdots ,r_{n-1} ), n \in \mathbb{Z}, k \in \mathbb{ N}$.
\end{definition}

\begin{theorem}\label{th:2.1}
For  $n \in \mathbb{Z}, k \in \mathbb{N}$, the numbers  $\check{C}h_{n;\bar{\pmb \alpha }, \bar{\pmb r}}^{(k)}$   satisfy the  following  relation
\begin{equation} \label{2.2}
\check{C}h_{n;\bar{\pmb \alpha}, \bar{\pmb r}}^{(k)}=\sum_{m=0}^{|r|}   s_{\bar{\pmb \alpha}} (n, m; \bar{\pmb r} ) \sum_{\ell_1=0}^m \cdots \sum_{\ell_k=0}^m \prod_{i=0}^k S(m,\ell_i )  Ch_{\ell_i }.
\end{equation}
\end{theorem}

\begin{proof}
From Eq. (\ref{2.1}), we have
\begin{eqnarray}\label{2.4}
\check{C}h_{n;\bar{\pmb \alpha}, \bar{\pmb r}}^{(k)}&=&
\int_{\mathbb{Z}_p}\cdots \int_{\mathbb{Z}_p} \sum_{m=0}^{|r|}  s_{\bar {\pmb \alpha}} (n,m;\bar{\pmb r}) (x_1 x_2\cdots x_k  )^m d\mu_{-1} (x_1 )\cdots d\mu_{-1} (x_k )
\nonumber\\
&=&\sum_{m=0}^{|r|}  s_{\bar{\pmb \alpha}} (n,m;\bar{\pmb r})\int_{\mathbb{Z}_p} \cdots \int_{\mathbb{Z}_p}(x_1 x_2\cdots x_k  )^m  d\mu_{-1} (x_1 )\cdots d\mu_{-1} (x_k )
\nonumber\\
&=&\sum_{m=0}^{|r|} s_{\bar{\pmb \alpha}}  (n,m;\bar{\pmb r})  \int_{\mathbb{Z}_p} (x_1  )^m  d\mu_{-1} (x_1 )\cdots \int_{\mathbb{Z}_p} (x_k  )^m  d\mu_{-1} (x_k )
\nonumber\\
 &=&\sum_{m=0}^{|r|} s_{\bar{\pmb \alpha}} (n,m;\bar{\pmb r} )\left[ \int_{\mathbb{Z}_p} \sum_{\ell_1=0}^m S(m,\ell_1 )  (x_1  )_{\ell_1}  d\mu_{-1} (x_1 )\cdots \int_{\mathbb{Z}_p} \sum_{\ell_k=0}^m S(m,\ell_k )  (x_k  )_{\ell_k }  d\mu_{-1} (x_k )\right]
 \nonumber\\
&=&\sum_{m=0}^{|r|}  s_{\bar{\pmb \alpha}}  (n,m;\bar{\pmb r} ) \left[ \sum_{\ell_1=0}^m S(m,\ell_1 )  \int_{\mathbb{Z}_p} (x_1  )_{\ell_1}  d\mu_{-1} (x_1 )\cdots \sum_{\ell_k=0}^m S(m,\ell_k )  \int_{\mathbb{Z}_p} (x_k  )_{\ell_k}  d\mu_{-1} (x_k ) \right]
\nonumber\\
&=&\sum_{m=0}^{|r|}  s_{\bar{\pmb \alpha}} (n,m;\bar{\pmb r})  \left[ \sum_{\ell_1=0}^m S(m,\ell_1 )  Ch_{\ell_1}\cdots \sum_{\ell_k=0}^m S(m,\ell_k )  Ch_{\ell_k } \right]
\nonumber\\
&=&\sum_{m=0}^{|r|}  s_{\bar{\pmb \alpha}}  (n,m;\bar{\pmb r})  \left[ \sum_{\ell_1=0}^m \sum_{\ell_2=0}^m \cdots \sum_{\ell_k=0}^m S(m,\ell_1 )S(m,\ell_2 )\cdots S(m,\ell_k ) Ch_{\ell_1} \cdots Ch_{\ell_k } \right],
\end{eqnarray}

\noindent
then Equation (\ref{2.2}) is obtained and this completes the proof.
\end{proof}

\begin{theorem} \label{th:2.2.}
The numbers $\check{C}h_{n;\bar{\pmb \alpha},\bar{\pmb r}}^{(k)}$   satisfy the  relation
\begin{equation} \label{2.5}
\check{C}h_{n;\bar{\pmb \alpha},\bar{\pmb r}}^{(k)}= \sum_{m=0}^{|r|}  s_{\bar{\pmb \alpha}}  (n,m;\bar{\pmb r})  \sum_{\ell_1=0}^m\cdots \sum_{\ell_k=0}^m \prod_{i=0}^k \frac{(-1)^{\ell_i} \ell_i ! S(m,\ell_i )}{\ell_i+1} ,
\end{equation}

\noindent
which gives a relationship of multiparameter Chaghee numbers of the first kind in terms of the multiparameter non central Stirling numbers of the second kind and Stirling number of the first kind, see \cite{8,9,10}.
\end{theorem}

\begin{proof} Substituting Eq. (\ref{2.1}) into (\ref{2.2}), we obtain (\ref{2.5}).
\end{proof}

\begin{remark} \label{rem:1}
\begin{eqnarray} \label{2.6}
\int_{\mathbb{Z}_p}  \int_{\mathbb{Z}_p} \cdots \int_{\mathbb{Z}_p} \prod_{i=0}^{n-1} (x_1 x_2\cdots x_k-\alpha_i  )^{r_i } d\mu_{-1} (x_1 ) \cdots d\mu_{-1} (x_k )
 \nonumber\\
=\sum_{m=0}^{|r|}   s_{\bar{\pmb \alpha}}  (n,m;\bar{\pmb r} ) \sum_{\ell_1=0}^m \cdots \sum_{\ell_k=0}^m \prod_{i=0}^k \frac{(-1)^{\ell_i } \ell_i! S(m,\ell_i )}{\ell_i+1} .
\end{eqnarray}
\end{remark}

\begin{definition} \label{def:2.1}
The multiparameter higher order Changhee polynomials   $\check{C}h_{n;\bar{\pmb \alpha}, \bar{\pmb r}}^{(k)}(x)$  are defined by
\begin{equation} \label{2.7}
\check{C}h_{n;\bar{\pmb \alpha},\bar{\pmb r}}^{(k)} (x)=\int_{\mathbb{Z}_p} \cdots \int_{\mathbb{Z}_p} \prod_{i=0}^{n-1} (x_1 x_2\cdots x_k x-\alpha_i  )^{r_i} d\mu_{-1} (x_1 )\cdots d\mu_{-1} (x_k ).
\end{equation}
\end{definition}

\subsection{Some special cases}
\noindent
{\bf Case 1:} (i) Setting  $r_i=r,\alpha_i=i$  in Eq. (\ref{2.7}), we have
\begin{eqnarray*}
\check{C}h_{n;i,\bar{\pmb r}}^{(k)} (x)=\int_{\mathbb{Z}_p} \cdots \int_{\mathbb{Z}_p} \prod_{i=0}^{n-1} (x_1 x_2\cdots x_k x-i )^r  d\mu_{-1} (x_1 )\cdots d\mu_{-1} (x_k )
\\
=\int_{\mathbb{Z}_p} \cdots \int_{\mathbb{Z}_p}(x_1 x_2 \cdots x_k x )_{nr} d\mu_{-1} (x_1 ) \cdots d\mu_{-1} (x_k ).
\end{eqnarray*}

\noindent
Replacing  $nr$ by $n$, we obtain the higher-order Changhee polynomials, which defined by Kim see \cite{7}.

\noindent
(ii) Setting  $r_i=r,\alpha_i=i$  in Eq. (\ref{2.1}), we obtain
\begin{equation*}
\check{C}h_{n;i,\bar{\pmb r}}^{(k)}=\int_{\mathbb{Z}_p} \cdots \int_{\mathbb{Z}_p} (x_1 x_2\cdots x_k  )_{nr}  d\mu_{-1} (x_1 )\cdots d\mu_{-1} (x_k ).
\end{equation*}

\noindent
Replacing  $nr$ by $n$ we obtain the higher-order Changhee numbers.\\

\noindent
{\bf Case 2:} (i) Setting  $r_i=r,\alpha_i=\alpha$ in Eq. (\ref{2.7}), we obtain
\begin{eqnarray*}
\check{C}h_{n;\alpha,r}^{(k)} (x) &=&
\int_{\mathbb{Z}_p} \cdots \int_{\mathbb{Z}_p} (x_1 x_2 \cdots x_k x-\alpha )^{nr} d\mu_{-1} (x_1 )\cdots d\mu_{-1} (x_k )
\\
&= & \int_{\mathbb{Z}_p} \cdots \int_{\mathbb{Z}_p} \sum_{\ell=0}^{nr} S(nr,\ell)  (x_1 x_2 \cdots x_k x-\alpha )_{\ell} d\mu_{-1} (x_1 )\cdots
d\mu_{-1} (x_k )
\\
&= & \sum_{\ell=0}^{nr} S(nr,\ell)  \int_{\mathbb{Z}_p} \cdots \int_{\mathbb{Z}_p} (x_1 x_2 \cdots x_k x-\alpha )_{\ell} d\mu_{-1} (x_1 )\cdots
d\mu_{-1} (x_k )
\\
&=& \sum_{\ell=0}^{nr} S(nr,\ell) \check{C}h_{\ell,\alpha }^{(k)} (x).
\end{eqnarray*}

\noindent
(ii) Setting  $r_i=r, \alpha_i=\alpha$ in Eq. (\ref{2.1}), we obtain
\begin{eqnarray}	
\check{C}h_{n;\alpha,r}^{(k)}&=&
\int_{\mathbb{Z}_p} \cdots \int_{\mathbb{Z}_p} (x_1 x_2 \cdots x_k-\alpha)^{nr} d\mu_{-1} (x_1 )\cdots d\mu_{-1} (x_k )
\nonumber\\ \label{2.8}
&=& \sum_{\ell=0}^{nr} S(nr,\ell) \check{C}h_{\ell,\alpha}^{(k)}.
\end{eqnarray}

\noindent
At $\alpha_i=0$ in  Eq. (\ref{2.8}), we obtain
\begin{equation*}
\check{C}h_{n;0,r}^{(k)} (x)= \sum_{\ell=0}^{nr} S(nr,\ell) \check{C}h_\ell^{(k) } (x).
\end{equation*}

\noindent
At $\alpha_i=0$ in  Eq. (\ref{2.8}), we obtain
\begin{equation*}
\check{C}h_{n;0,r}^{(k)} = \sum_{\ell=0}^{nr} S(nr,\ell) \check{C}h_\ell^{(k)}.
\end{equation*}

\noindent
{\bf Case 3:} (i) Setting  $r_i=1, \alpha_i=\alpha$ in Eq. (\ref{2.7}), we obtain
\begin{equation*}
\check{C}h_{n;\alpha,1}^{(k)} (x) =\sum_{\ell=0}^n S(n,\ell) \check{C}h_{\ell,\alpha}^{(k)} (x).
\end{equation*}

\noindent
(ii) Setting  $r_i=1, \alpha_i=\alpha$ in Eq. (\ref{2.1}), we obtain
\begin{equation*}
\check{C}h_{n;\alpha,1}^{(k)}=\sum_{\ell=0}^n S(n,\ell) \check{C}h_{\ell,\alpha}^{(k) }.
\end{equation*}

\noindent
(iii) Setting  $r_i=1,  \alpha_i=1$ in Eq. (\ref{2.7}), we obtain
\begin{eqnarray*}
\check{C}h_{n;1,1}^{(k)} (x) & = & \int_{\mathbb{Z}_p} \cdots \int_{\mathbb{Z}_p} (x_1 x_2\cdots x_k x-1 )^n  d\mu_{-1} (x_1 )\cdots d\mu_{-1} (x_k )
\\
&= & \sum_{\ell=0}^n S(n,\ell) \check{C}h_{\ell,1}^{(k)} (x).
\end{eqnarray*}

\noindent
(iv) Setting  $r_i=1,   \alpha_i=1$ in Eq. (\ref{2.1}), we obtain
\begin{eqnarray*}
\check{C}h_{n;1,1}^{(k)} &=&
\int_{\mathbb{Z}_p} \cdots \int_{\mathbb{Z}_p} (x_1 x_2\cdots x_k-1 )^n  d\mu_{-1} (x_1 )\cdots d\mu_{-1} (x_k )
\\
&= & \sum_{\ell=0}^n S(n,\ell) \check{C}h_{\ell,1 }^{(k)}.
\end{eqnarray*}

\noindent
{\bf Case 4:} (i) Setting  $r_i=1,\alpha_i=0$ in Eq. (\ref{2.7}), we obtain
\begin{eqnarray*}
Ch_{n;0,1}^{(k) } (x) & = &
\int_{\mathbb{Z}_p} \cdots \int_{\mathbb{Z}_p} (x_1 x_2 \cdots x_k x)^n d\mu_{-1} (x_1 )\cdots d\mu_{-1} (x_k )
\\
&=& \sum_{\ell=0}^n S(n,\ell) Ch_{\ell}^{(k)} (x).
\end{eqnarray*}

\noindent
(ii) Setting  $r_i=1, \alpha_i=0$ in Eq. (\ref{2.1}), we obtain
\begin{equation*}
Ch_{n;0,1}^{(k)}= \sum_{\ell=0}^n S(n,\ell) Ch_{\ell}^{(k)}.
\end{equation*}

\noindent
{\bf Case 5:} (i) Setting  $r_i=1, \alpha_i=i,i=0,1,\cdots,n-1$  in Eq. (\ref{2.7}), we have
\begin{eqnarray}
Ch_{n;i,1}^{(k) } (x) & = &
\int_{\mathbb{Z}_p} \cdots \int_{\mathbb{Z}_p} \prod_{i=0}^{n-1} (x_1 x_2\cdots x_k x-i )  d\mu_{-1} (x_1 )\cdots d\mu_{-1} (x_k )
\nonumber\\
&= &
\int_{\mathbb{Z}_p} \cdots \int_{\mathbb{Z}_p} (x_1 x_2 \cdots x_k x)_n  d\mu_{-1} (x_1 )\cdots d\mu_{-1} (x_k )
\nonumber\\ \label{2.10}
&=& Ch_n^{(k) } (x),
\end{eqnarray}

\noindent
we obtain the higher-order Changhee polynomials which defined by Kim see \cite{7}.

\noindent
(ii) Setting  $r_i=1, \alpha_i=i,i=0,1,\cdots,n-1$  in Eq. (\ref{2.1}), we have
\begin{eqnarray}
Ch_{n;i,1}^{(k)} & = &
\int_{\mathbb{Z}_p} \cdots \int_{\mathbb{Z}_p} \prod_{i=0}^{n-1} (x_1 x_2 \cdots x_k-i ) d\mu_{-1} (x_1 )\cdots d\mu_{-1} (x_k )
\nonumber\\
&=& \int_{\mathbb{Z}_p} \cdots \int_{\mathbb{Z}_p} (x_1 x_2 \cdots x_k )_n  d\mu_{-1} (x_1 )\cdots d\mu_{-1} (x_k )
\nonumber\\ \label{2.11}
& = & Ch_n^{(k)},
\end{eqnarray}

\noindent
we obtain the Changhee numbers of the first kind with order $k$, see \cite{7}.\\

\noindent
{\bf Case 6:}  Setting  $x_1 x_2\cdots x_k=x$  in Eq. (\ref{2.1}), we obtain
\begin{equation} \label{2.12}
Ch_{n;\bar{\pmb \alpha},\bar{\pmb r}} =\int_{\mathbb{Z}_p} (x-\alpha_0  )^{r_0} (x-\alpha_1  )^{r_1} \cdots  (x-\alpha_{n-1}  )^{r_{n-1}} d\mu_{-1} (x).
\end{equation}

\begin{corollary}\label{cor:1}
The numbers $Ch_{n;\bar{\pmb \alpha}, \bar{\pmb r}}$  satisfy the  relation
\begin{equation} \label{2.13}
Ch_{n;\bar{\pmb \alpha}, \bar{\pmb r}} = \sum_{i=0}^{|r|} S(n,i; \bar{\pmb \alpha}, \bar{\pmb r})  Ch_i .
\end{equation}
\end{corollary}

\begin{proof}
Let   $x_1 x_2\cdots x_k=x$  in Eq. (\ref{2.2}), we obtain Eq. (\ref{2.13}).
\end{proof}

\begin{corollary} \label;{cor:2}
The numbers $Ch_{n;\bar{\pmb \alpha}, \bar{\pmb r}}$  satisfy the  relation
\begin{equation} \label{2.14}
Ch_{n;\bar{\pmb \alpha}, \bar{\pmb r}} = \sum_{\ell=0}^{|r|} s_{\bar{\pmb \alpha}}  (n,\ell;\bar{\pmb r})  E_{\ell}.
\end{equation}
\end{corollary}

\begin{proof}
From Eq. (\ref{1.12}) and Eq. (\ref{2.13}), we obtain Eq. (\ref{2.14}).
\end{proof}

\noindent
{\bf Case 7:} Setting  $r_i=1$ in Eq. (\ref{2.12}), we obtain
\begin{equation} \label{2.15}
Ch_{n;\bar{\pmb \alpha}} =\int_{\mathbb{Z}_p}(x-\alpha_0  )(x-\alpha_1  )\cdots (x-\alpha_{n-1}  )  d\mu_{-1} (x),
\end{equation}

\noindent
which we define $Ch_{n;\bar{\pmb \alpha}}$ by generalized Changhee numbers of the first kind.

\begin{corollary}\label{cor:3}
The numbers $Ch_{n;\bar{\pmb \alpha}}$  satisfy the  relation
\begin{equation} \label{2.16}
Ch_{n;\bar{\pmb \alpha}} = \sum_{i=0}^n S(n,i;\bar{\pmb \alpha} )   Ch_i.
\end{equation}
\end{corollary}

\begin{proof}
Let   $r_i=1,i=0,1,\cdots,n-1$  in Eq. (\ref{2.13}), we obtain Eq. (\ref{2.16}).
\end{proof}

\noindent
{\bf Case 8:} Setting $\int_{\mathbb{Z}_p} \cdots \int_{\mathbb{Z}_p} =\int_0^{\ell_1} \cdots \int_0^{\ell_k}$ in (\ref{2.1}), we obtain
\begin{equation} \label{2.17}
C_{n;\bar{\pmb \alpha}, \bar{\pmb r}}^{(k) }=\int_0^{\ell_1} \cdots \int_0^{\ell_k} \prod_{i=0}^{n-1} (x_1 x_2 \cdots x_k-\alpha_i  )^{r_i } dx_1\cdots dx_k,
\end{equation}

\noindent
multiparameter poly-Cauchy numbers of the first kind.

\section{Multiparameter Chaghee Numbers of the Second Kind}
\begin{definition}\label{def:3.1}
The multiparameter Chaghee numbers of the second kind  $\widehat{Ch}_{n;\bar{\pmb \alpha},\bar{\pmb r}}^{(k)}$ are defined by
\begin{equation} \label{3.1}
\widehat{Ch}_{n;\bar{\pmb \alpha},\bar{\pmb r}}^{(k)}= \int_{\mathbb{Z}_p} \cdots \int_{\mathbb{Z}_p} \prod_{i=0}^{n-1} (-x_1 x_2\cdots x_k-\alpha_i )^{r_i } d\mu_{-1} (x_1 ) \cdots d\mu_{-1} (x_k ).
\end{equation}
\end{definition}

\begin{theorem} \label{th:3.1}
The numbers $\widehat{Ch}_{n;\bar{\pmb \alpha},\bar{\pmb r}}^{(k)}$ satisfy the  relation
\begin{equation} \label{3.2}
\widehat{Ch}_{n;\bar{\pmb \alpha},\bar{\pmb r}}^{(k)}= \sum_{m=0}^{|r|}  s_{\bar{\pmb \alpha}}  (n,m;\bar{\pmb r})  \sum_{\ell=0}^m L(m,n)  \sum_{\ell_1=0}^m \cdots \sum_{\ell_k=0}^m \prod_{i=0}^k S(m,\ell_i ) Ch_{\ell_i },
\end{equation}

\noindent
where $L(m,n)$ is the Lah numbers, see \cite{13}.

\end{theorem}

\begin{proof}
From Eq. (\ref{3.1}) we have
\begin{eqnarray}
\widehat{Ch}_{n;\bar{\pmb \alpha}, \bar{\pmb r}}^{(k)} &= &
\int_{\mathbb{Z}_p}\cdots \int_{\mathbb{Z}_p} \sum_{m=0}^{|r|}  s_{\bar{\pmb \alpha}}  (n,m;\bar{\pmb r}) (–x_1 x_2 \cdots x_k  )_m d\mu_{-1} (x_1 )\cdots d\mu_{-1} (x_k )
\nonumber\\ \label{3.4}
&=&
\int_{\mathbb{Z}_p} \cdots \int_{\mathbb{Z}_p} \sum_{m=0}^{|r|}  s_{\bar{\pmb \alpha}} (n,m;\bar{\pmb r}) \sum_{\ell=0}^m L(m,n) (x_1 x_2\cdots x_k  )_{\ell} d\mu_{-1} (x_1 )\cdots d\mu_{-1} (x_k ).
\end{eqnarray}

\noindent
Substituting Eq. (\ref{2.2}) into (\ref{3.4}), we obtain (\ref{3.2}.
\end{proof}

\begin{definition}\label{def:3.1}
The multiparameter higher order Changhee polynomials   $\widehat{\check{C}h}_{n;\bar{\pmb \alpha}, \bar{\pmb r}}^{(k)}$ are defined by
\begin{equation} \label{3.5}
\widehat{\check{C}h}_{n;\bar{\pmb \alpha},\bar{\pmb r}}^{(k)} (x)=\int_{\mathbb{Z}_p}\cdots \int_{\mathbb{Z}_p} \prod_{i=0}^{n-1} (-x_1 x_2\cdots x_k x-\alpha_i  )^{r_i} d\mu_{-1} (x_1 )\cdots d\mu_{-1} (x_k ).
\end{equation}
\end{definition}

\subsection{Some special cases}
{\bf Case 1:} (i) Setting  $r_i=r,\alpha_i=i$  in Eq. (\ref{3.5}), we have
\begin{eqnarray*}
\widehat{\check{C}h}_{n;i,\bar{\pmb r}}^{(k)} (x) & =&
\int_{\mathbb{Z}_p} \cdots \int_{\mathbb{Z}_p} \prod_{i=0}^{n-1} (-x_1 x_2\cdots x_k x-i )^r d\mu_{-1} (x_1 )\cdots d\mu_{-1} (x_k )
\\
&=& \int_{\mathbb{Z}_p} \cdots \int_{\mathbb{Z}_p} (-x_1 x_2\cdots x_k x )_{nr} d\mu_{-1} (x_1 )\cdots d\mu_{-1} (x_k ).
\end{eqnarray*}

\noindent
Replacing  $nr$ by $n$, we obtain the higher-order Changhee polynomials, which defined by Kim see \cite{7}.

\noindent
(ii) Setting  $r_i=r,\alpha_i=i$  in Eq. (\ref{3.1}), we obtain
\begin{equation*}
\widehat{\check{C}h}_{n;i,\bar{\pmb r}}^{(k)}= \int_{\mathbb{Z}_p} \cdots \int_{\mathbb{Z}_p} (-x_1 x_2\cdots x_k  )_{nr} d\mu_{-1} (x_1 ) \cdots d\mu_{-1} (x_k ).
\end{equation*}

\noindent
Replacing  $nr$ by $n$, we obtain the higher-order Changhee numbers.\\

\noindent
{\bf Case 2:} (i) Setting  $r_i=r,\alpha_i=\alpha$ in Eq. (\ref{3.5}), we obtain	
\begin{eqnarray}
\widehat{\check{C}h}_{n;\alpha,r}^{(k)} (x) &= &
\int_{\mathbb{Z}_p} \cdots \int_{\mathbb{Z}_p} (-x_1 x_2\cdots x_k x-\alpha )^{nr} d\mu_{-1} (x_1 )\cdots d\mu_{-1} (x_k )
\nonumber\\
&=& \int_{\mathbb{Z}_p} \cdots \int_{\mathbb{Z}_p} \sum_{\ell=0}^{nr} S(nr,\ell)  (-x_1 x_2\cdots x_k x-\alpha )_{\ell}  d\mu_{-1} (x_1 )\cdots d\mu_{-1} (x_k )
\nonumber\\
&=& \sum_{\ell=0}^{nr} S(nr,\ell) \int_{\mathbb{Z}_p} \cdots \int_{\mathbb{Z}_p} (-x_1 x_2 \cdots x_k x-\alpha )_{\ell}  d\mu_{-1} (x_1 )\cdots d\mu_{-1} (x_k )
\nonumber\\ \label{3.6}
&=& \sum_{\ell=0}^{nr} S(nr,\ell) \widehat{\check{C}h}_{\ell,\alpha} )^{(k) } (x).
\end{eqnarray}

\noindent
(ii) Setting  $r_i=r,\alpha_i=\alpha$ in Eq. (\ref{3.1}), we obtain
\begin{eqnarray}
\widehat{\check{C}h}_{n;\alpha,r}^{(k)} &= &
\int_{\mathbb{Z}_p} \cdots \int_{\mathbb{Z}_p} (-x_1 x_2\cdots x_k-\alpha )^{nr} d\mu_{-1} (x_1 )\cdots d\mu_{-1} (x_k ),
\nonumber\\   \label{3.7}
&=& \sum_{\ell=0}^{nr} S(nr,\ell) \widehat{\check{C}h}_{\ell,\alpha}^{(k) }.
\end{eqnarray}

\noindent
At $\alpha_i=0,r_i=r$ in  Eq. (\ref{3.5}), we obtain
\begin{equation*}
\widehat{\check{C}h}_{n;0,r}^{(k)} (x)=\sum_{\ell=0}^{nr}S(nr,\ell) \widehat{\check{C}h}_{\ell}^{(k)} (x).
\end{equation*}

\noindent
At $\alpha_i=0,r_i=r$ in  Eq. (\ref{3.2}), we obtain
\begin{equation*}
\widehat{\check{C}h}_{n;0,r}^{(k)}=\sum_{\ell=0}^{nr} S(nr,\ell)  \widehat{\check{C}h}_{\ell}^{(k)}.
\end{equation*}

\noindent
{\bf Case 3:} (i) Setting  $r_i=1,\alpha_i=\alpha$ in Eq. (\ref{3.6}), we obtain
\begin{equation*}
\widehat{\check{C}h}_{n;\alpha,1}^{(k)} (x)=\sum_{\ell=0}^n S(n,\ell) \widehat{\check{C}h}_{\ell,\alpha }^{(k)} (x).
\end{equation*}

\noindent
(ii) Setting  $r_i=1,\alpha_i=\alpha$ in Eq. (\ref{2.2}), we obtain
\begin{equation*}
\widehat{\check{C}h}_{n;\alpha,1}^{(k) }=\sum_{\ell=0}^n S(n,\ell) \widehat{\check{C}h}_{\ell,\alpha }^{(k)}.
\end{equation*}

\noindent
(iii) Setting  $r_i=1,  \alpha_i=1$ in Eq. (\ref{3.7}), we obtain
\begin{eqnarray*}
\widehat{\check{C}h}_{n;1,1}^{(k)} (x) &= & \int_{\mathbb{Z}_p} \cdots \int_{\mathbb{Z}_p} (-x_1 x_2 \cdots x_k x-1 )^n  d\mu_{-1} (x_1 )\cdots d\mu_{-1} (x_k )
\\
&= & \sum_{\ell=0}^n S(n,\ell) \widehat{\check{C}h}_{\ell,1}^{(k) } (x).
\end{eqnarray*}

\noindent
(iv) Setting  $r_i=1,  \alpha_i=1$ in Eq. (\ref{2.2}), we obtain
\begin{eqnarray*}
\widehat{\check{C}h}_{n;1,1}^{(k)} &= & \int_{\mathbb{Z}_p} \cdots \int_{\mathbb{Z}_p} (-x_1 x_2\cdots x_k-1 )^n  d\mu_{-1} (x_1 )\cdots d\mu_{-1} (x_k )
\\
&=& \sum_{\ell=0}^n S(n,\ell) \widehat{\check{C}h}_{\ell,1}^{(k)}.
\end{eqnarray*}

\noindent
{\bf Case 4:} (i) Setting  $r_i=1,\alpha_i=0$ in Eq. (\ref{3.6}), we obtain
\begin{eqnarray*}
\widehat{\check{C}h}_{n;0,1}^{(k)} (x) & =&
\int_{\mathbb{Z}_p} \cdots \int_{\mathbb{Z}_p} (x_1 x_2 \cdots x_k x)^n  d\mu_{-1} (x_1 )\cdots d\mu_{-1} (x_k )
\\
&=& \sum_{\ell=0}^n S(n,\ell) \widehat{\check{C}h}_{\ell}^{(k)} (x).
\end{eqnarray*}

\noindent
(ii) Setting  $r_i=1,\alpha_i=0$ in Eq. (\ref{2.2}), we obtain
\begin{equation*}
\widehat{\check{C}h}_{n;0,1}^{(k) }= \sum_{\ell=0}^n S(n,\ell) \widehat{\check{C}h}_{\ell}^{(k)}.
\end{equation*}

\noindent
{\bf Case 5:} (i) Setting  $r_i=1, \alpha_i=i,i=0,1,\cdots,n-1$  in Eq. (\ref{3.6}), we have
\begin{eqnarray}
\widehat{\check{C}h}_{n;i,1}^{(k) } (x) & = &
\int_{\mathbb{Z}_p} \cdots \int_{\mathbb{Z}_p} \prod_{\ell=0}^{n-1} (-x_1 x_2 \cdots x_k x-i )   d\mu_{-1} (x_1 )\cdots d\mu_{-1} (x_k )
\nonumber\\
& = & \int_{\mathbb{Z}_p} \cdots \int_{\mathbb{Z}_p} (-x_1 x_2\cdots x_k x)_n  d\mu_{-1} (x_1 )\cdots d\mu_{-1} (x_k )
\nonumber\\ \label{3.8}
& =& \widehat{\check{C}h}_n^{(k) } (x),
\end{eqnarray}

\noindent
we obtain the higher-order Changhee polynomials, which defined by Kim see \cite{7}.

\noindent
(ii) Setting  $r_i=1, \alpha_i=i,i=0,1,\cdots, n-1$  in Eq. (\ref{3.1}), we have
\begin{eqnarray*}
\widehat{\check{C}h}_{n;i,1}^{(k) } &=&
\int_{\mathbb{Z}_p} \cdots \int_{\mathbb{Z}_p} \prod_{i=0}^{n-1} (-x_1 x_2 \cdots x_k-i )    d\mu_{-1} (x_1 )\cdots d\mu_{-1} (x_k )
\\
&=& \int_{\mathbb{Z}_p} \cdots \int_{\mathbb{Z}_p} (-x_1 x_2 \cdots x_k )_n  d\mu_{-1} (x_1 )\cdots d\mu_{-1} (x_k )
\\
&=& \widehat{\check{C}h}_n^{(k) },
\end{eqnarray*}

\noindent
we obtain the Changhee numbers of the first kind with order $k$, see \cite{7}.\\

\noindent
{\bf Case 6:} Setting  $-x_1 x_2 \cdots x_k=-x$  in Eq. (\ref{3.1}), we obtained
\begin{equation} \label{3.9}
\widehat{\check{C}h}_{n;\bar{\pmb \alpha},\bar{\pmb r}} = \int_{\mathbb{Z}_p} (-x-\alpha_0  )^{r_0} (-x-\alpha_1  )^{r_1} \cdots (-x-\alpha_{n-1} )^{r_{n-1}} d\mu_{-1} (x).
\end{equation}

\begin{corollary} \label{cor:61}
The numbers $\widehat{\check{C}h}_{n;\bar{\pmb \alpha}, \bar{\pmb r}}$ satisfy the  relation

\begin{equation} \label{3.10}
\widehat{\check{C}h}_{n;\bar{\pmb \alpha},\bar{\pmb r}} = \sum_{i=0}^{|r|} (-1)^i S_{\bar{\pmb \alpha}}  (n,i;\bar{\pmb r} )  Ch_i.
\end{equation}
\end{corollary}

\noindent{\bf Case 7:} Setting  $r_i=1$ in Eq. (\ref{3.9}), we obtain
\begin{equation*}
\widehat{\check{C}h}_{n;\bar{\pmb \alpha}}= \int_{\mathbb{Z}_p} (-x-\alpha_0  )(-x-\alpha_1) \cdots (-x-\alpha_{n-1} ) d\mu_{-1} (x).
\end{equation*}

\begin{corollary}\label{cor:3}
The numbers $\widehat{\check{C}h}_{n;\bar{\pmb \alpha}}$ satisfy the  relation
\begin{equation} \label{3.11}
\widehat{\check{C}h}_{n;\bar{\pmb \alpha}} = \sum_{i=0}^n (-1)^i S_{\bar{\pmb \alpha}}  (n,i)  Ch_i.
\end{equation}
\end{corollary}

\noindent
{\bf Case 8:} Setting   $\int_{\mathbb{Z}_p} \int_{\mathbb{Z}_p}\cdots \int_{\mathbb{Z}_p}=\int_0^{\ell_1} \int_0^{\ell_2} \cdots \int_0^{\ell_k}  $ in (\ref{3.1}), we obtain
\begin{equation} \label{3.12}
\widehat{C}_{n;\bar{\pmb \alpha}, \bar{\pmb r}}^{(k)} = \int_0^{\ell_1}  \int_0^{\ell_2} \cdots \int_0^{\ell_k} \prod_{i=0}^{n-1} (-x_1 x_2 \cdots x_k-\alpha_i  )^{r_i} dx_1 dx_2 \cdots dx_k,
\end{equation}

\noindent
multiparameter poly-Cauchy numbers of the second kind.


\end{document}